\documentclass[12pt]{amsart}
\newtheorem{thm}{Theorem}[section]
\newtheorem{prop}[thm]{Proposition}

\newtheorem{cor}[thm]{Corollary}

\newtheorem{defn}[thm]{Definition}

\newtheorem{exam}[thm]{Example}
\begin{document}


\newcommand{\lb}{\left \langle }
\newcommand{\rb}{\right \rangle}
\newcommand{\bi}{\bibitem}
\newcommand{\ci}{\cite}
\newcommand{\no}{\nonumber}
\newcommand{\be}{\begin{equation}}
\newcommand{\ee}{\end{equation}}
\newcommand{\vsp}{\vspace{.1in}}
\newcommand{\hsp}{\hspace{.2in}}
\newcommand{\bea}{\begin{eqnarray}}
\newcommand{\eea}{\end{eqnarray}}
\newcommand{\bean}{\begin{eqnarray*}}
\newcommand{\eean}{\end{eqnarray*}}
\newcommand{\ba}{\begin{array}}
\newcommand{\ea}{\end{array}}
\newcommand{\ol}{\overline}
\newcommand{\mt}{{ \mathcal T }}

\newcommand{\xgn}{{ X_{g,n}}}
\newcommand{\xg}{{ X_{g}}}
\newcommand{\dzgs}{{{\bf \mathcal  Z ^{g*} }}}
\newcommand{\zgs}{{{\bf Z^{g*}}}}
\newcommand{\dzg}{{{\bf \mathcal  Z ^g }}}
\newcommand{\dz}{{{\bf \mathcal  Z  }}}
\newcommand{\dza}{{{\bf \mathcal  Z ^a }}}
\newcommand{\dzb}{{{\bf \mathcal  Z ^b }}}

\newcommand{\zg}{{{\bf Z^g}}}
\newcommand{\zas}{{{\bf Z^{a*}}}}
\newcommand{\zbs}{{{\bf Z^{b^*}}}}
\newcommand{\za}{{{\bf Z^a}}}
\newcommand{\zb}{{{\bf Z^b}}}
\newcommand{\zag}{{{\bf Z^g_{ab}}}}
\newcommand{\zah}{{{\bf Z^h_{ab}}}}
\newcommand{\lab}{{L^2_{ab}}}
\newcommand{\Lab}{{L^2{[0,\frac{1}{b}]}\times L^2{[0,a]}}}

\newcommand{\ja}{{ \frac{-j}{a} }}
\newcommand{\na}{{ \frac{-n}{a} }}
\newcommand{\fkb}{{\frac{k}{b}}}
\newcommand{\fjb}{{\frac{j}{b}}}
\newcommand{\ka}{{ \frac{-k}{a} }}
\newcommand{\fab}{{ \left(\frac{a}{b}\right)}}
\newcommand{\sfb} {{\sqrt{\frac{1}{b}} }}
\newcommand{\fb}{{ \frac{1}{b}}} 
\newcommand{\lr}{{ L^2 (\mathbb R)}}
\newcommand{\la}{{ L^2[0,a]}}
\newcommand{\mr}{{\mathbb R}}
\newcommand{\mz}{{ \mathbb Z}}
\newcommand{\mc}{{ \mathbb C}}
\newcommand{\hc}{{ L^{\infty}_1\left(\ell_2\right)}}
\newcommand{\LB}{{ L^{\infty}_{\fb}\left(\ell_2\right)}}

\newcommand{\cp}{{\mathcal{P}}}
\newcommand{\bd}{\begin{displaystyle}}
\newcommand{\ed}{\end{displaystyle}}
\newcommand{\g}{{g_{m,n}}}

   \renewcommand{\theequation}{\arabic{section}.\arabic{equation}}
\def\ldots{\mathinner{\ldotp\ldotp\ldotp}}
\def\ldots{\mathinner{\cdotp\cdotp\cdotp}} \def \norm{| \! | \! |}
\def \N{\rm {\bf N}} \def \R{\rm {\bf R}} \def \K{\rm {\bf K}} \def
\diag{\text{diag }} \def \tr{\text{tr }} \def \Com{\text{Com }} \def
\cal{\mathcal} \def \Bbb{\mathbb} \def \det{\text{det }}

\title{ A  Hilbert $C^*$-Module for  Gabor systems }  
\author{Michael  Coco and  M.C. Lammers}
\email{Michael Coco, University of South Carolina, coco@math.sc.edu}
\email{M.C. Lammers, University of South Carolina,
lammers@math.sc.edu}
\maketitle
\begin{abstract} We construct  Hilbert $C^*$-modules useful for
studying Gabor systems and show that they are  Banach Algebra under
pointwise multiplication.  For rational $ab<1$ we prove that the set of
functions $g \in \lr$ so that $(g,a,b)$ is a Bessel system is an ideal
for the Hilbert $C^*$-module given this pointwise algebraic structure.
This allows us to give a multiplicative  perturbation theorem for frames. 
Finally we show that a system $(g,a,b)$ yields a frame for $\lr$ iff
it is a modular frame for the given Hilbert $C^*$-module.
  
 \end{abstract}

\section{Intro}  A {\bf Gabor system} for a function $ g \in \lr$
 is obtained by applying modulations and translations ($g(x-ka)e^{2\pi
 i m b x}$) to this function.  In the study of these systems many
 results have been obtained by using the function-valued inner product
 $$\lb f,g\rb_1(x)=\sum_{k \in\mz} f(x-k)\ol{g(x-k)}.$$ Many of the
 results from Hilbert space theory can be reproduced in this setting.
 However, one of the most basic properties of a Hilbert space is lost.
 Namely, for $f,g,h\in \lr$ the function $\lb f,g\rb_1(x) h(x) $ need
 not be be in $\lr$.  This brings us to the study of {\bf Hilbert
 $C^*$-Modules (HCM)}.  A Hilbert $C^*$-Module is the generalization
 of the Hilbert space in which the inner product maps into a
 $C^*$-algebra.  By restricting to a subspace of $\lr$ we produce a
 Hilbert $C^*$-Module for the inner product above containing the
 functions $g$ which are of most interest in the study of Gabor
 systems.\\

Our study of this HCM and its relation to Gabor systems is organized
in the following manner.  In section~2 we provide some of the known
results involving Gabor systems and Gabor frames along with some of
the basics regarding HCM's.  In section~3 we define our HCM's and
develop a list of basic properties which relate to Gabor systems. In
section 4 we show that our HCM is a Banach algebra under pointwise
multiplication and show that the set of functions with a finite upper
frame bound form an ideal for this algebra. In addition we give a
Balian-Low type theorem in the positive direction and show that one
can apply an additive perturbation result of Christenson and Heil to
obtain a multiplicative one in our Banach algebra.  In section~5 we
make the connection between {\bf modular frames} for our HCM (defined
by Frank and Larson) and Gabor frames for $\lr$.  We show that the
system $(g,a,b)$ is a frame for $\lr$ iff the translates of $g$,
$\{T_{ka}g\}_{k\in \mz}$, form a modular frame for our HCM.  \\

The authors would like to thank Michael Frank for his many helpful
suggestions and references regarding Hilbert $C^*$-modules.
Throughout all sums are considered to be over $\mz$ or $\mz^@$
 unless otherwise stated.


\section{Preliminaries}

In 1952 Duffin and Schaeffer \cite{DS} defined frames:

\begin{defn}  A sequence $(f_{n})_{n\in \Bbb Z}$ of elements
of a Hilbert space $H$ is called a {\bf frame} if there are constants
$A,B > 0$ such that
\begin{equation}
A\|f\|^{2}\le \sum_{n\in \Bbb Z}|\lb f,f_{n}\rb |^{2}\le B\|f\|^{2},\
\
\text{for all}\ \ f\in H. \no
\end{equation}
\end{defn}

For a function $f \in \lr$ and $a,b\in \mr$ define the operators:

 \[
\text{Translation: }  T_{a}f(x) = f(x-a), \qquad 
\text{ Modulation: }   E_{b}f(x) = e^{2{\pi}ibx}f(x), \,
\]

\begin{defn}
If $a,b\in \mr$ and $g\in \lr$ we let $(E_{mb}T_{na}g)=\g$ and call
$\left \{(E_{mb}T_{na}g)\right \}_{m,n\in \mz}$ a {\bf Gabor system}
(also called a Weyl-Heisenberg system) and denote it by $(g,a,b)$. If
this system is a frame then we call it a {\bf Gabor frame}.  We denote
by $(g,a)$ the family $(T_{na}g)_{n\in \mz}$.
\end{defn}

Now we wish to define some important operators associated with a Gabor
system $(g,a,b)$. Let $(e_{m,n})$ be an orthonormal basis for
$\lr$. We call the operator $\mathcal T_g:\lr\rightarrow \lr$ given by
$\mathcal T_ge_{m,n}=E_{mb}T_{na}g$ the {\bf preframe operator}
associated with $(g,a,b)$.  The adjoint, $\mathcal T_g^*$, is called
the {\bf frame transform} and for $f\in \lr$ and $m,n\in \Bbb Z$ we
have $\lb
\mathcal T_g^{*}f,e_{m,n}\rb = \lb f,\mathcal T_g(e_{m,n})\rb =\lb
f,(E_{mb}T_{na}g)\rb $.  Thus
\[  
\mathcal T_g^{*}f = \sum_{m,n}\lb f,E_{mb}T_{na}g\rb e_{m,n} \,\text {
and } \, \|\mathcal T_g^{*}f\|^{2} = \sum_{m,n}|\lb f,E_{mb}T_{na}g\rb
|^{2} \text{ for all } f\in \lr. \] It follows that the preframe
operator is bounded if and only if $(g_{m,n})$ has a finite upper
frame bound $B$.  Finally we define the {\bf frame operator associated
with $(g,a,b)$}:
\[ S(f)=\mathcal T_g\mathcal T_g^*(f)=\sum_{m,n\in \mz}\lb
f,E_{mb}T_{na}g\rb (E_{mb}T_{na}g).\] It is known that if $(g,a,b)$ is
a frame then $S$ is a positive, invertible operator.  We will be
interested in when there is a finite upper frame bound for a Gabor
system.

\begin{defn} Let $\{\psi_n\}_{n \in \mz}$ be an orthonormal basis for
a Hilbert space $H$.  A {\bf Bessel sequence} is any sequence
$\{\phi_n\}_{n \in\mz}$(i.e. not necessarily one generated by the
operators $E_{mb}$ and $T_{na}$) so that $T(\sigma)=\sum_{n \in
\mz}\lb
\sigma,\psi_n\rb_H  \phi_n$ defines a bounded operator from $H$ to
$H$.  For fixed $a$ and $b$ we denote the set of functions $g$ so that
$\{E_{mb}T_{na}g\}$ is a Bessel sequence by $X_Z$.\\
\end{defn}

We will use a bracket product notation compatible with Gabor systems.
We define the {\bf $a$-inner product} for $f,g, \in \lr$ and $a\in
\mr^+$:
\[\lb f,g\rb _a=\lb f,g\rb _a(x)=\sum_k f(x-ka)\ol{g(x-ka)}\text{ and }
 \|f\|_a(x)=\sqrt{\lb f,f\rb _a(x)}. \]

For a detailed development of this $a$-inner product which includes two
forms of a Riesz Representation we refer the reader to \cite{CL} but
here we give a brief list of some of the properties which will be
useful later.

\begin{prop}\label{pip} For $f,g \in \lr $ \\
(a) $\lb f,f\rb _a\ge 0$ and = 0 iff $f=0$ \\ (b) $\lb f,g\rb _a=
\overline{\lb g,f\rb _a} $ \\ (c) $\lb \phi f,g\rb _a= \phi \lb f,g\rb
_a$ for all periodic functions such that $f,\phi f \in \lr$\\ (d) $\lb
f,h+g\rb _a = \lb f,h\rb_a +\lb f,g\rb_a$\\ (e) $|\lb f,g\rb_a|^2 \le
\|f\|_a \|g\|_a $\\ (f) $\lb f,g\rb _a \in L^1[0,a]$ \\ (g) $\lb
f,g\rb _a$ is a-periodic on $\mr$\\ (h) $\|f\|^2_\lr =\int_0^a \lb
f,f\rb_a dx$\\ (i) $\|f+g\|_a \le \|f\|_a +\|g\|_a$
\end{prop}

While we see that the $a$-inner product has many of the same
properties as the usual inner product, there is one major difference
as we mentioned in the introduction.  Consider $f,g,h \in\lr$ where \[
f(x)=g(x)=
\frac{1}{x^\frac{1}{3}}{\bf 1}_{[0,1]},\,\, h(x)={\bf 1}_{[0,1]}.\] Then
$\lb f,g \rb_1 (x)h(x)= \bd \frac{1}{x^\frac{2}{3}}\ed {\bf 1}_{[0,1]}$
and this function is clearly not in $\lr$.\\

We say that an operator $L:\lr \to \lr$ is {\bf $a$-factorable} if
 $L(\phi f)=\phi L(f)$ for all bounded $a$-periodic functions $\phi$. The operators associated with Gabor systems are
$\fb$-factorable and hence have the following representations which we
refer to as {\bf compressions} since we compress the modulation into
the $\fb$-inner product.

\begin{prop}\cite{CL,rs2} If $\|g\|_{\fb}\le B$ a.e., then the frame
operator and the preframe operator have the following compressions
\[\mathcal T_g(f)= \sfb \sum_k\lb f,e_k\rb _{\fb}T_{ka}g, \quad
\mathcal T_g^*(f)= \sfb \sum_k\lb f,T_{ka} g\rb _{\fb}e_k \text { and
}\]
\[S_g(f)=\fb \sum_k\lb f,T_{ka}g\rb_\fb T_{ka}g\]
where $e_k =T_{\fkb} 1_{[0,\fb)}$ .
\end{prop}

\begin{defn}For $\lambda > 0$ the {\bf Zak transform} of a function 
$f\in\lr$ is \[{ Z_\lambda}(f)(t,v)= \lambda^{1/2}\sum_{k \in \mz} f(
\lambda(t-k))e^{2 \pi i k v}, \, t,v \in \mr.\] We may interpret this
as a unitary operator from $\lr$ to $L^2[0,1]^2$.  Since we are
frequently interested in the case where $\lambda=1$, we will also use
the notation $Z_1=Z$.
\end{defn}

Before we present some of the basics regarding Hilbert $C^*$-modules
we give two very important results for $(g,1,1)$ systems. The first is
due to Daubechies and the second was proved independently by Balian and
Low.  A nice treatment of both can be found in \cite{FS}.

\begin{thm}\label{DJ}\cite{D} The system $(g,1,1)$ yields  a frame with 
frame bounds A and B iff
\[A \le |Z(f)(t,v)| \le B \text{ for a.e. } (t,v)\in [0,1]^2\]
\end{thm}

\begin{thm}\label{BL}{\bf(Balian-Low)} If the system $(g,1,1)$
yields a frame, then

\[ \|tg(t)\|_\lr \, \|\gamma\hat g(\gamma)\|_\lr=+\infty \]
\end{thm}

Now we give some of the basics regarding Hilbert $C^*$-modules.  We
refer the reader to \cite{rw} for a thorough treatment of this
subject.  Essentially, a Hilbert $C^*$-module is the generalization of
a Hilbert space in which the inner product takes values in a
$C^*$-algebra instead of the complex numbers.

\begin{defn}  Let $A$ be a $C^*$-algebra.  An {\bf inner product
$A$-module } is an $A$-module $\mathcal M$ with a mapping $\lb
\cdot,\cdot
\rb_A:\mathcal M\times \mathcal M\to A$ so that for all $x,y,z \in
\mathcal M$ and  $a\in A$ \\

a) $\lb x,x \rb_A\ge 0 $ (as an element of $A$);\\ b) $\lb x,x \rb_A
=0$ iff $x=0$;\\ c) $\lb x,y \rb_A = \lb y,x \rb_A^*$; \\ d) $\lb ax,x
\rb_A= a \lb x,x \rb_A$ ; \\ e) $\lb x+y,x \rb_A =\lb x,z \rb_A +\lb
y,z \rb_A.$\\ 

If the space $\mathcal M$ is complete with respect to
the norm $\|x\|^2_A=\|
\lb x,x \rb_A\|_A $, then we say $\mathcal M$ is a {\bf Hilbert $C^*$-module
with respect to $A$} or more simply a {\bf Hilbert $A$-module}. When
referring to a Hilbert $C^*$-module we will also use the abbreviation
{\bf HCM}.
\end{defn}

We provide a very simple example.

\begin{exam} Let $A$ be a $C^*$-algebra.  The $A$-valued inner product
defined by $\lb a,b \rb_A=ab^*$ makes $A$ a Hilbert $A$-module.
\end{exam}


\section{The HCM}

Given the way we have defined the $a$-inner product and the Hilbert
$C^*$-module it is natural to try to use the $a$-inner product to turn
$\lr$ into a HCM.  Unfortunately, the $a$-inner product maps $\lr$
into $L^1[0,a]$ which is not a $C^*$-algebra.  For this reason we
consider the subspace of $\lr$ consisting of those functions which are
mapped into $L^\infty [0,a]$.  This yields the following Lebesgue-Bochner space.\\

For $b>0$ let $\LB$ be defined to be the set of measurable functions
$f:\mr \to \mc$ for which the norm
\[\|f\|^2_{\LB} =esssup_{[0,\fb]} \sum_k|f(x-\fkb )|^2 \]is finite.
 Proposition
\ref{pip} gives us that $\LB$ is a Hilbert $C^*$-Module ({\bf HCM})
over the $C^*$-algebra $L^\infty [0,\fb)$ with $C^*$-valued inner
product and $C^*$-valued norm
\[\lb f,g\rb_\fb(x)= \sum_kf(x-\fkb )\ol{g(x-\fkb )}
\,\text{ and } \, \|f\|_\fb(x)=\left(\sum_k|f(x-\fkb )|^2\right)^{1/2} .\]

We begin with a few basic remarks:\\

{\bf 1}.  Weak-* convergence. Because we are constructing our $HCM$
with Gabor systems in mind we will only require the $C^*$-valued inner
product to converge weak-* in $L^\infty[0,1]$.  We motivate this by
considering the case$a=b=1$.  In this case we have the classification
Theorem \ref{DJ}.  With this theorem it is easy to see that the system
$(g,1,1)$ with $\bd g=\sum_{k=0}^\infty {\bf1}_{[k+\frac{1}{2^{k+1}},
k+\frac{1}{2^k})} \ed$ produces a Gabor frame because $|Z(g)(t,v)|=1$
everywhere. However the series $\lb g,g
\rb_1$ does not converge in norm in $L^\infty[0,1]$. If we only
require that the series converge in the weak* topology we avoid this
problem.  This assumption results in frames being what Frank and
Larson \cite{FL} have termed nonstandard modular frames.  We will
address this more in section 5.\\

{\bf 2}. Proposition~\ref{pip}~(h) yields the following inequality
which shows that $\LB$ embeds continuously in $\lr$
\[\|f\|^2 =\int _0^\fb |\lb f,f\rb_\fb| dx 
\le \sup_{[0,1]} |\lb f,f\rb_\fb| =\|f\|_\LB^2.\]
     
Finally, since $\LB$ contains continuous functions of compact support,
we see it is norm dense in $\lr$.\\

 {\bf 3}. $\LB\ne \lr\cap L^\infty(\mr)$.  We illustrate this for
 $\hc$ and note that it is easily generalized.  Let $f_k= {\bf
 1}_{[k,k+\frac{1}{k^2})}$ and $f=\sum_{k=0}^1 f_k$.  Then $f \in
 \lr\cap L^\infty(\mr)$ but $f\not \in \hc$ for $\lb f,f\rb_1 (x)$ is
 unbounded near zero. \\

{\bf 4}. We interpret a basic result about Gabor frames in this space.
A necessary condition for $\{\g \}_{m,n\in \mz}$ to be a Gabor frame
in $\lr$ with upper frame bound $B$ is that $esssup_{[0,\fb]}
\sum_k |g(x-\fkb )|^2 \le B $ a.e..  In other words, if $\g$ is a frame
for $\lr$ with upper frame bound $B$ then $g$ in the $B$-ball of
$\LB$.\\

{\bf 5}.  The Wiener amalgam spaces $W(L,G)$ are defined using the
norm of two spaces where the local behavior is determined by (L) and
the global behavior by (G).  For example,

\[W(L^p,\ell^q)=\left\{ f : \|f\|_{W(L^p,\ell^q)}
=\left ( \sum_k\|f\cdot e_k\|^q_p \right )^\frac{1}{q} <\infty \right
\}.\]

We will refer to ${\bf W}(\mr)=W(L^\infty,\ell^1)$ as the Wiener
Algebra.  The Segal algebra may be viewed as a specific Wiener amalgam
space:

\[{\bf S}_0(\mr)=W(C_0,\ell^1)=\left\{ f \in W(L^\infty,\ell^1) 
 : f \text{ is continuous } \right \}.\]

These algebras have been used effectively in the study of Gabor
systems by Benedetto, Heil, Walnut, Feichtinger, Zimmerman and
Strohmer ( See chapters 2,3 and 8 of \cite{FS}). From the inequality
\[ \|g\|_\hc \le \|g\|_{W(L^\infty, \ell_2)} \le  \|g\|_{W(L^\infty,
\ell_1)}\] 
we conclude that the Segal algebra and the Weiner algebra embed in
$\hc$.  That is, for any $g$ in ${\bf S}_0(\mr)$ we have

\[ \|g \|_\hc \le  \|g\|_{{\bf W}(\mr)}\le \|g\|_{{\bf S}_0(\mr)}\]

{\bf 6.}  One may view the Zak transform here as the analogue of the
Fourier transform when $a=b=1$.  We give the usual notation for this
but add the representation with respect to the inner product for our
HCM.  We should point out that although these transforms agree for
$t,v \in [0,1]^2$, the latter is periodic in $t$ and $v$ while the
former yields only a quasi periodicity.

\[Z(f)(t,v)= \sum_k  f(t+k)e^{-2 \pi i k v}
= \sum_k \lb f,e_k \rb_1(t)e^{-2 \pi i k v}\] where $e_k=T_k{\bf
1}_{[0,1)}.$

\[\hat f(v)=\int_\mr f(t)e^{-2 \pi i t v}dt =\int_0^1 Z(f)(t,v)e^{-2
\pi i t v}dt\] 

 {\bf 7.}  The Fourier transform is not a bounded linear operator on
 $\hc$. Indeed, if we let $f(x)={\bf 1}_{[0,1]}
\bd \frac{1}{x^\frac{1}{3}}\ed$ it can be shown that the absolute
value of the Fourier transform is bounded by the function
\[ g(x)=  \begin{cases}  2  &x\in [-1,1]\\
\frac{2}{x^\frac{2}{3}}  & \text{otherwise}   \end{cases} . \]

The function $g(x)$ is clearly in $\hc$ which implies $\hat f$ is
also.  However, since the original function $f$ is not even bounded,
it is clearly not in $\hc$.  Thus by the properties of the
Fourier transform we see that $\hat{ \hat f} $ is not in $\hc$.

         
\section{Bessel Systems}

Before we examine the role Bessel systems play in our HCM's we show
that these HCM's are Banach algebras under pointwise multiplication.

\begin{prop} The space $\LB$ is a Banach algebra under pointwise
multiplication.

\end{prop}

\begin{proof}
 The justification can be made with one inequality.  Note that the
 inequality below is not derived by Cauchy-Shwarz, but by the fact
 that all the terms are positive and the left hand side is just the
 diagonal of the product of the sums on the right.
\[\sum_k|f  g(x-\fkb)|^2 \le \sum_k|f(x-\fkb)|^2 \sum_k|g(x-\fkb)|^2 .\]
\end{proof}

  Unfortunately $\LB$ is not a $C^*$-algebra. It is also clear that
  this Banach Algebra has no identity since $f(x)={\bf 1}_\mr $ is
  clearly not in $\LB$.  \\
  
Now we turn our attention to the functions that yield Bessel systems.
We show that the functions which produce a finite upper frame bound
for $a=b=1$ form a Banach space which may be embedded in
$L^\infty_1(\ell_2)$. For now let us refer to this as the Zak space,
denoted $X_Z$.  More formally we define the space,

\[X_Z=\left\{f| f\in \lr\, \text{and}\, \|f\|_{X_Z}= esssup_{t,v \in [0,1]^2
}|Z(f)(t,v)|<\infty\right \}.\]

Since $Z$ is a unitary operator from $\lr$ to $L^2[0,1]^2$ it is easy
to see that $\|\cdot\|_{X_Z}$ is actually a norm.  The essential
inequality in the following proposition is well known and is usually
stated in the form $G_0=\sum_k|g(x-k)|^2 <B$ .  We prove
it here for completeness and for the development of the Zak transform
on $\hc$.

\begin{prop} \label{Zembed}The Zak space embeds continuously  
in $\hc$, that is
\[\|g\|_\hc \le \|g\|_{X_Z} \text{ for all } g\in X_Z. \] \end{prop}

\begin{proof}  To see that $X_Z$ embeds in $\hc$  we look at what the Zak
transform does on $\hc$.  Let us define

\[L^\infty
L^2[0,1]^2=\left\{F(t,v)|\sup_{t\in[0,1]}\left(\int_0^1|F(t,v)|^2dv\right)
^\frac{1}{2} <\infty\right \}.\] Now let $f\in \hc$ so

\bean \|Z(f)(t,v)\|_{L^\infty L^2[0,1]^2}&=&\sup_{t\in[0,1]}\left(\int_0^1
|\sum\lb f,e_k\rb_1 e^{-2 \pi i k v}|^2dv\right)^\frac{1}{2}\\
&=&\sup_{t\in[0,1]}
\left(\sum|\lb f,e_k\rb_1|^2\right)^\frac{1}{2}=\|f\|_\hc
\eean It is easy to see from here that the Zak transform is an
isometry between these two spaces. Further, since the $L^2[0,1]$ norm
is bounded by the $L^\infty[0,1]$ norm, we see that $\|g\|_\hc \le
\|g\|_{X_Z}$ for $g\in X_Z$.

\end{proof}

Our next example shows that $X_Z \ne \hc$.

\begin{exam}There exists $g\in \hc $ with $g\not \in X_Z$ \end{exam}

\begin{proof}
 Consider the function $g=\bd \sum_{n=1}^\infty \frac{e_n}{n}\ed$. Clearly
 $g\in \hc $. Now Corollary 3.7 from \cite{CCJ} states that a positive
 real-valued function is in $X_Z$ iff $\bd \sum_k|\lb g, T_k g \rb_1| \le
 \infty\ed$.  However computation shows that $\bd \lb g, T_k g
 \rb_1=(1/k)\sum_{n=1}^{k}\frac{1}{n}\ed$ for the above $g$. These are
 square summable but not summable. \end{proof}

As we mentioned a detailed study of the Balian and Low theorem
(Theorem \ref{BL}) can be found in \cite{FS}.  In this study an
amalgam version of the BLT theorem for the Segal algebra ${\bf S}_0$
due to Heil is given.

\begin {thm}\cite{H}  Let $g\in \lr$. If $(g,1,1)$  forms a frame
for $\lr$, then
\[g\not \in {\bf S_0} \text{ and } \hat g \not \in {\bf S_0}\]

\end{thm}

A direct Corollary of the continuous embedding yields the positive
Balian-Low type result which in light of remark {\bf 7} from the
previous section is non-trivial.

\begin{cor} If the system $(g,1,1)$ is a frame then $g\in \hc$ and
$\hat g \in \hc$
\end{cor}
 
The proof follows directly from the theorem above and the fact that
$(g,1,1)$ forms a frame if and only if $(\hat g,1,1)$ does.  This,
however, does not characterize the space $X_Z$.  We showed in the
proof of Proposition
\ref{Zembed} that $Z$ is an isometry from  $\hc $ to $L^\infty
L^2[0,1]^2$.  Recall that $Z(\hat g)=e^{2
\pi i t v}Z(g)(v,-t)$.  It suffices to show that there exists
a function $F(t,v)\in L^\infty L^2[0,1]^2$ so that $F(v,-t) \in
L^\infty L^2[0,1]^2 $ yet $F(t,v) \not \in L^\infty[0,1]^2$.  The
function $F(t,v)=\bd v^{-\frac{t}{3}}\ed$ is such a function. Given
that the Segal algebra (see remark {\bf 5} section 3) embeds in $\hc$
we may interpret these results as follows.  If $(g,1,1)$ forms a frame
then $g,\hat g\in \hc\setminus{\bf S_0}$. \\

Now we show that the functions with finite upper frame bound have an algebraic connection to
$L_1^\infty(\ell_2)$. We recall an algebraic term.  If $M$ is a ring without
identity we say the {\bf square ideal}, $M^2$ is the ideal of the ring
formed by taking finite sums of products from $M$.  That is
\[M^2=\left\{\sum_{i=1}^n x_iy_i: x_i,y_i \in M\right\}.\]

The theorem below not only shows that the Zak space is an ideal in
$\hc$ but it contains the square ideal.

\begin{thm} For any $f,g\in L_1^\infty(\ell_2)$ their product $fg$
has a finite upper frame bound.  In particular, the Zak space is an
ideal for $L_1^\infty(\ell_2)$.
\end{thm}

\begin{proof} Let $f,g$ in $L^\infty(\ell_2)$ and consider $Z(fg)$.

\bea| Z(fg)(t,v)|^2&=&\left |\sum_k \lb fg,e_k\rb_1(t)e^{-2 \pi i k
v}\right|^2\no \\ &\le & \sum_k|f(t+k)|^2\sum_k|g(t+k)|^2 \le
\|f\|_{\hc}
\|g\|_{\hc}. \no \eea
\end{proof}

As usual the argument may be mimicked in cases where $ab$ is rational.
To get the most general rational case one probably needs to use the
Zak matrices of Zibulski-Zeevi which may be found in
\cite{FS} section 1.5.  We present the case where $a=\frac{1}{2},b=1 $
and leave the rest to the reader.

\begin{thm} If $f,g$ in $L_1^\infty(\ell_2)$, then $(fg,\bd
\frac{1}{2}\ed ,1)$
is a Bessel system.
\end{thm}

\begin{proof} To see this we apply a standard technique when using
the  Zak transform. Applying the Zak transform to the 
frame operator of the system $(g,\frac{1}{2},1)$  we get:   

\bean  Z(S_g)(f)&=&Z(\sum_k \lb f,T_{\frac{k}{2}} g \rb_1 T_{\frac{k}{2}}g) \\
 &=& Z(\sum_k \lb f,T_{{k}} g \rb_1 T_{k}g +
 Z(\sum_k \lb f,T_{k} T_{\frac{1}{2}}g \rb_1 T_{k}(T_{\frac{1}{2}}g)\\
 &=& Z(\sum_k \lb f,T_{{k}} g \rb_1 T_{k}g )+Z(
 \sum_k \lb f,T_{k} T_{\frac{1}{2}}g \rb_1 T_{k}(T_{\frac{1}{2}}g)\\
 &=& Z(f)\left(|Z(g)|^2+|Z(T_{\frac{1}{2}}g)|^2\right).
\eean

 This implies that the  system $(g,\frac{1}{2},1)$
has a finite upper frame bound iff
$|Z(g)|^2+|Z(T_{\frac{1}{2}}g)|^2 <B$. In view of this and the
proposition above all we need to show is that $|Z(T_{\frac{1}{2}}fg)|$
bounded.  This follows easily from the fact that
$\|T_{\frac{1}{2}}f\|_\hc =\|f\|_\hc$
\end{proof}

One may interpret this theorem another way.  We state it as a corollary
in the case $a=b=1$.

\begin{cor} For any $g\in L_1^\infty(\ell_2)$ and $x\in \mr$ 
 the operator
\[\mathcal V^1_g(f)(t,v,x)  =Z(f T_xg)(t,v) e^{-2\pi itv}\]
 is a bounded operator from $\LB$ to $L^\infty (\mr^3)$.
\end{cor}

\begin{proof} First we point out that because we are considering the periodic
extension of the Zak transform in $t,v\in[0,1]^2$ we really only need
$[0,1]^2\times \mr$ in place of $\mr^3$.  Given $x\in
\mr$ we have $\|T_xg\|_\LB=\|g\|_\LB$.  The result follows from the above inequality.
\end{proof}

Let us point out that in many ways $\mathcal{V}^1_g(f)(t,v,x)$ resembles the
short time Fourier transform
\[\mathcal V_g(f,x)=\int_\mr f(t)\ol {g(t-x)}e^{-2 \pi i tv}dt\]
For this reason we term this operator the {\bf windowed Zak
transform}. 

We have a simple corollary for producing frames in $\hc$. Again we
state it here for the case $a=b=1$ but it is easily generalized to
many rational cases.

\begin{cor} If $f\in \hc$ and $\inf |Z(f^2)| >A$ then $(f^2,1,1)$ is a
frame. \end{cor}

Finally we end the section by applying the following perturbation
theorem of Christensen and Heil to produce frames.

\begin{thm}\cite{CH} Let $H$ be a Hilbert space, $\{f_i\}_{i=1}^\infty $ be
a frame with bounds $A,B$ and let $\{g_i\}_{i=1}^\infty \subset H$. If
there exists $R<A$ so that  
\[ \sum_{i=1}^\infty|\lb h,f_i-g_i\rb_{\mathcal H}|^2 \le 
R \, \|h\|^2 \text{ for all } h \in H,\] then
$\{g_i\}_{i=1}^\infty $ is a frame with bounds
$A(1-\sqrt{\frac{R}{A}})^2$ and $B(1-\sqrt{\frac{R}{A}})^2$.
\end{thm}

In the Gabor case this means that if $(g,a,b)$ is a frame then so is
$(f,a,b)$ if the system $(f-g,a,b)$ has an upper frame bound less than
that of $(g,a,b)$.  We use this theorem and one of the inequalities
above to produce a multiplicative perturbation result.

\begin{prop} Let $ab<1$ and $ab$ rational. If $(g,a,b)$ produces a
frame with bounds $A,B$ and if $(f-1)\in \hc $ with $\|f-1\|^2_\hc \le
\frac{R}{B}$, then $(fg,a,b)$ is a frame.
\end{prop}

\begin{proof}  As usual, to highlight the essential components 
 of the argument we prove the case $a=b=1$. Let us point out that
 $f-1$ need not be in $\hc$ even if $f$ is, since $1\notin \hc$. 
 By the theorem above it
 is enough to show that $fg-g$ has a finite upper frame bound less
 than $A$.  In view of Theorem \ref{DJ} it is enough to show
 $|Z(fg-g)|<R$.

\bea| Z(fg-g)(t,v)|^2&=&\left |\sum_k\lb(f-1)g,e_k \rb_1(t)e^{-2 \pi i k
v}\right|^2\no \\ &\le & \sum_k|(f-1)(t+k)|^2\sum_k|g(t+k)|^2 \no\\
&\le& \|f-1\|^2_{\hc}
\|g\|^2_{\hc} \le \frac{R}{B} \cdot B .\no \eea

Where the fact that $\|g\|_{\hc} \le B$ follows from Proposition
\ref{Zembed}.  In this case we get frame bounds
$A(1-\sqrt{\frac{R}{A}})^2$ and $B(1-\sqrt{\frac{R}{A}})^2$.

\end{proof}  


\section{$a$-frames and modular frames}

An immediate concern is whether the operators $\mathcal T_g $ and
$\mathcal T^*_g $ are bounded operators on $\LB$ if they are bounded
operators on $\lr$.  The answer is yes and follows from Proposition
5.9 of \cite{CL} which we state below.

\begin{prop}If $L:\lr \to \lr $ is an a-factorable operator, then $L$
is bounded iff for all $f \in \lr $
\[\|L(f)\|_a(t) \le \|L\|\|f\|_a(t)\] where
$\|L\|=\sup_{\|f\|_\lr=1}\|L(f)\|_\lr $.
\end{prop}

Given the compressed representations of $\mathcal T_g $ and $\mathcal
T^*_g $ it is easy to see that they are $\fb$-factorable on $\lr$.
Since the inequality in the proposition holds for all $f\in \lr$ it
certainly holds for $f\in \LB$ and we have our result by taking
supremums on both sides.\\

We now introduce two definitions both of which are abstractions of
frames.  The first is that of an $a$-frame for $\lr$ and the second is
the analogue of a frame in a Hilbert $C^*$-module.\\
 
\begin{defn}\cite{CL}  A sequence $\{f_n\} \in \lr$ is an {\bf $a$-frame for
$\lr$ } if there exists $A,B$ such that for all $f\in\lr$
\bea \label{a1} A \|f\|^2_a(x) \le \sum_n |\lb f,f_n \rb_a(x)|^2 \le
B\|f\|_a^2(x) \text{ a.e..}\eea
\end{defn} 

Since we are dealing with frames for different spaces we will add the
term ``modular'' to Frank and Larson's definition of a frame for a
HCM.

\begin{defn} \cite{FL} Let $\mathcal M$ be a HCM over the  unital
$C^*$-algebra $A$.  A sequence $\{x_n\} \in \mathcal M$ is said to be
a {\bf modular frame } for $\mathcal M$ if there are real constants
$C,D>0$ such that
\[C\lb x,x\rb_A \le \sum_n \lb x,x_n\rb_A\lb x_n,x\rb_A \le D\lb
x,x\rb_A\] for all $x\in \mathcal M$. If the sum in the middle of the
inequality always converges in norm this is referred to as a {\bf
standard modular frame} and if the middle sum converges weakly for
some $x\in
\mathcal M$ we will call this a {\bf non-standard modular frame}. 
\end{defn}

We point out here that our results are in a different direction than
the ``standard'' results given by Frank and Larson.  Since we only
require the sum $$\sum_k f(x-\fkb)\ol {g(x-\fkb)}$$ to converge weakly
in $L^\infty[0,\fb]$, we are dealing with the nonstandard case.  Most
of their examination dealt with standard modular frames in HCM's which
were at worst countably generated.  Neither of these conditions is met
in our case. Now we give the result of Casazza and Lammers regarding
the connection between $a$-frames and Gabor frames.

\begin{thm}\cite{CL} Let 
$g_n(x)=g(x-na)$.  Then $\{\g \}_{m,n\in \mz}$ is a Gabor frame for
$\lr$ with frame bounds $A,B$ iff $\{g_n\}$ is a $\fb$-frame for
$\lr$. \end{thm}

This is somewhat surprising because this says that the frame
inequality holds pointwise for the bracket product.  We use this
result to show the connections between Gabor frames and modular
frames.

\begin{thm}For $g\in \lr$,  $\{\g\}_{m,n \in \mz}$ is  a Gabor frame
for $\lr$ iff $\{g_n\}_n$ is a non-standard modular frame for $\LB$.
\end{thm}

\begin{proof} $\Rightarrow$ This follows directly from the theorem
above since $\LB$ is a subspace of $\lr$.\\

 $\Leftarrow$ For any $f\in \lr$ consider $f_0 =\bd
\frac{f}{\|f\|_{\fb}} \ed$. Then $\lb f_0,f_0\rb_\fb \le 1$ and
hence $f_0 \in \LB$. Since $\{g_n\}$ is a modular frame for $\LB$,

\bea A \|f_0\|^2_\fb(x) \le \sum_n |\lb f_0,g_n \rb_\fb(x)|^2 \le
B\| f_0\|_\fb^2(x) \text{ a.e.} \no \\ A
\frac{\|f\|^2_\fb(x)}{\|f\|^2_\fb(x)} \le \frac{1}{\|f\|^2_\fb(x)}
\sum_n |\lb f,g_n \rb_\fb(x)|^2 \le B
\frac{\|f\|^2_\fb(x)}{\|f\|^2_\fb(x)} \text{ a.e.} \no \\ A
\|f\|^2_\fb(x) \le \sum_n |\lb f,g_n \rb_\fb(x)|^2 \le B\|
f\|_\fb^2(x) \text{ a.e..} \no\eea

So we are done by the previous theorem.
\end{proof}

Finally we mention that there is a convolution for our HCM in the case
$a=b=1$.  This is developed in much greater detail in
\cite{ML}. Our convolution is simply the preframe operator
$\mt_g(f)$.  That is, we define $g*_1f=\mt_g(f)$ because
\[  Z(\mt_g(f))=\sum_k \lb f, e_k\rb_1 Z(T_kg)=\sum_k \lb f, e_k\rb_1
e^{2 \pi i k v}Z(g)=Z(g)Z(f).\] The reason we call this a convolution
is because the Zak transform turns it into multiplication. Clearly
this does not necessarily map to an element in $L^2[0,1]^2$ but from Holder's
inequality we see that $Z(\mt_g(f))\in L^1[0,1]^2$.



\end{document}